\documentclass[12pt]{article}
\input epsf      
\usepackage{graphicx}
\usepackage{amssymb}
\usepackage{amsfonts}
\usepackage{amsthm}
\usepackage{amsfonts}

\newcommand{\complex}{\mathbb{C}}

\newcommand{\half}{\mathbb{H}}

\newcommand{\compose}{\circ}

\newcommand{\ba}[1]{\begin{array}{#1}}
\newcommand{\ea}{\end{array}}

\newcommand{\be}{\begin{equation}}
\newcommand{\ee}{\end{equation}}
\newcommand{\bea}{\begin{eqnarray}}
\newcommand{\eea}{\end{eqnarray}}
\newcommand{\beann}{\begin{eqnarray*}}
\newcommand{\eeann}{\end{eqnarray*}}
\newcommand{\blength}{{b}}
\newcommand{\nterms}{{n}}

\def\reff#1{(\ref{#1})}
\begin{document}

\bibstyle{ams}

\title{
Computing the Loewner driving process of random curves in the half plane
}

\author{Tom Kennedy
\\Department of Mathematics
\\University of Arizona
\\Tucson, AZ 85721
\\ email: tgk@math.arizona.edu
}
\maketitle

\begin{abstract}

We simulate several models of random curves in the half plane and numerically
compute the stochastic driving processes that produce the curves 
through the Loewner equation.
Our models include models whose scaling limit is the 
Schramm-Loewner evolution (SLE) and models for which it is not.
We study several tests of whether the driving process is Brownian 
motion, as it is for SLE.
We find that testing only the normality of  the process at a 
fixed time is not effective at determining if the 
random curves are an SLE.
Tests that involve the independence of the increments of 
Brownian motion are much more effective.
We also study the zipper algorithm for numerically computing the driving
function of a simple curve. We give an implementation of this algorithm
which runs in a time $O(N^{1.35})$ rather than the usual $O(N^2)$,
where $N$ is the number of points on the curve. 
\end{abstract}

\bigskip

\newpage

\section{Introduction}
\label{intro}

The Loewner equation provides a means for encoding curves in the 
upper half plane that do not intersect themselves by a real-valued 
function. Let $\gamma(t)$ be such a simple curve with $0 \le t < \infty$. 
Let $\half$ denote the upper half of the complex plane, and let 
$\gamma[0,t]$ denote the image of $\gamma$ up to time $t$. Then 
$\half \setminus \gamma[0,t]$ is a simply connected domain. So there is 
a conformal map $g_t$ from this domain to $\half$. If the curve is 
suitably parametrized and $g_t$ is suitably normalized, 
then $g_t$ satisfies the differential equation 
\be
{\partial g_t(z) \over \partial t} = {2 \over g_t(z) - U_t},
\qquad g_0(z)=z
\ee
for some real valued function $U_t$ on $[0,\infty)$. 
The function $U_t$ is often called the driving function.

If our simple curve in the half plane is random, 
then the driving function $U_t$ is a stochastic process.
Schramm discovered that if the scaling limit of a two-dimensional 
model is conformally invariant and satisfies a certain Markov property, 
then this stochastic driving process must be a Brownian motion 
with mean zero \cite{schramm}. 
The only thing that is not determined is the variance. 
Schramm named this process stochastic Loewner evolution or SLE; it is now
often referred to as Schramm-Loewner evolution. 

Many critical two-dimensional models from statistical mechanics 
and probability satisfy these properties or are believed to satisfy them,
and so should be SLE for some $\kappa$ (the parameter that determines 
the variance). These include the loop-erased 
random walk \cite{lsw_lerw,zhan}, the self-avoiding walk \cite{lsw_saw}, 
interfaces in the critical Ising model \cite{smirnov_ising},
the Gaussian free field \cite{ss}, 
critical percolation \cite{smirnov_perc, cn}, 
and uniform spanning trees \cite{lsw_lerw}. 
More recent work has considered whether other models have random curves 
that are described by SLE.
The possibility that domain walls in spin glass ground states are SLE curves 
was studied numerically 
both by Amoruso, Hartman, Hastings, and Moore  \cite{ahhm} 
and by Bernard, Le Doussal, and Middleton \cite{bdm}.
Bernard, Boffetta, Celani and Falkovich 
considered simulations of certain isolines in two-dimensional turbulence 
\cite{bbcfA} and surface quasi-geostrophic turbulence \cite{bbcfB}.

In this paper we consider models which are definitely not SLE. 
They are based on well known lattice models - 
the loop-erased random walk (LERW), the self-avoiding walk (SAW) and 
the critical percolation exploration process. 
We distort these models by shrinking the random curves slightly 
in the vertical direction but not in the horizontal direction. 
In other words we apply a non-conformal transformation to the curves. 
Without distortion these models are all proven or conjectured to have 
a scaling limit given by SLE. 

One way to test if a model of random curves is SLE is to compute its
stochastic driving process and see if it is Brownian motion. 
In this paper we simulate these distorted models, 
numerically compute their stochastic driving process, and 
then test if they are Brownian motions. We also do this for the models
without distortion.
Our goal is to see how well one can determine whether or not a model 
is SLE by studying this stochastic driving process and to compare
various methods for testing if the stochastic driving process is a 
Brownian motion.

Another goal of this paper is to study the algorithm for computing 
the driving function of a given curve. The standard implementation of 
the ``zipper algorithm'' for doing this requires a time $O(N^2)$ where
$N$ is the number of points on the curve. We present an 
implementation that runs in a time $O(N^p)$ with $p$ approximately 
$1.35$. This implementation uses the same idea used in 
\cite{tk_sle} to simulate SLE curves quickly. 
We also study the difference in the driving function found using 
``tilted slits'' versus ``vertical slits'' and the effect of the 
number of points used on the curve to compute the driving function. 

\section{Distorted models}
\label{distorted}

We study the stochastic driving function of three models. We refer
to them as distorted models. For $\lambda>0$ we 
define a non-conformal map on the upper half plane by 
$\phi(x,y)=(x,\lambda y)$.
Given a model that produces random curves $\gamma$ in the upper half
plane, we consider the random curves $\phi \compose \gamma$. 
In other words, we stretch the curve by a factor of $\lambda$ in the 
vertical direction, but do not stretch it in the horizontal direction. 

We apply this distortion to the loop-erased random walk, the 
self-avoiding walk and percolation interfaces. For all three models
we consider the chordal version of the model in which the random curve
lies in the upper half plane and goes from the origin to $\infty$. 
For the LERW this means we take a half plane excursion and loop erase it. 
For the SAW this means we use the uniform probability measure on nearest 
neighbor walks with a fixed number of steps which begin at the origin and 
lie in the upper half plane. 
For percolation we consider site percolation on the triangular lattice 
in the upper half plane with boundary conditions which force an interface 
to start at the origin. Details of the definitions of these models and 
the parameters used in the simulations may be found in the appendix. 

There is no simple relation between the driving function for 
the curve $\gamma$ and the driving function for the distorted curve
$\phi \compose \gamma$. We study the driving function for the distorted model
as follows. We generate samples of the LERW, SAW or percolation interface
and then apply the distortion map $\phi$. Then we compute the driving
function of the distorted curve. The result is a collection of samples of the 
stochastic driving process of the distorted model. We then do various
statistical tests to see if this process is a Brownian motion.

We denote the driving process by $U_t$. 
All the models are invariant under reflections about the vertical 
axis. Hence $E[U_t]=0$. 
We begin by plotting the variance $E[U_t^2]$ as a function of $t$. 
We should emphasize that in the scaling limit,
all the models have a scaling property 
which implies that $E[U_t^2]$ is a linear function of $t$, even if 
the scaling limit is not an $SLE$.  
So this does not test whether $U_t$ is a Brownian motion. It only 
provides an estimate of $\kappa$ where $\kappa$ is the slope of the 
function $t \rightarrow E[U_t^2]$. 

The first statistical test is to see if the distribution of an individual
$U_t$ is normal. We use the Kolmogorov-Smirnov test. 
This test is based on the fact that for a continuous random variable $Y$,
if $F$ is the cumulative distribution of $Y$, then $F(Y)$ is uniformly 
distributed on $[0,1]$. Let $Y_1,Y_2, \cdots, Y_N$ be $N$ observations of 
the random variable $Y$, and let $Y_{(1)} < Y_{(2)} < \cdots Y_{(N)}$ 
be these numbers arranged in increasing order. Then the statistic is 
\be
D=\max_{ 1 \le k \le N} |F(Y_{(k)})-{k+{1 \over 2} \over N}| + { 1 \over 2N}
\ee
(This formula may appear different from that found in most statistics texts,
but it is in fact the same.) 
Under the null hypothesis that $Y$ has the distribution $F$, the limiting 
distribution of $\sqrt{N} D$ as $N \rightarrow \infty$ is known. 
For example, $P(\sqrt{N}D >1.36)$ is approximately $5 \%$. 
So if we compute this statistic for an individual $U_t$ with $F$ equal 
to the cumulative distribution for a normal random variable with mean 
zero and variance $\kappa t$, and find that the value of $\sqrt{N} D$ is large 
(say larger than $1.36$), then we conclude that the distribution of $U_t$ 
is not this normal distribution. We perform this Kolmogorov-Smirnov 
test for two values of the time, $T$ and $T/2$. 
The value of $T$ as well as other parameters used in the simulations
are given in the appendix.
The results are shown in the various tables in the columns 
labeled $D(T/2)$ and $D(T)$.

Our next test involves the independence of increments of Brownian motion. 
Let
\be
X_1=U_{T/2}, \quad X_2=U_T-U_{T/2}
\ee
If $U_t$ is a Brownian motion, then $X_1$ and $X_2$ are independent 
and have mean zero. So their product $X_1 X_2$ has mean zero. 
We test the hypothesis that its mean is zero. 
If $U_t$ is a Brownian motion, then the variance of $X_1 X_2$ can 
be found. Denote it by $\sigma^2$. The statistic we use, $Z$, is simply
the sample mean for $X_1 X_2$ divided by $\sigma/\sqrt{N}$. 
If $N$ is large and $U_t$ is a Brownian motion, then the distribution
of $Z$ is close to that of the standard normal. 
This test is shown in the tables in the column labeled $Z$.

For our last three tests we let $0 < t_1 < t_2 < \cdots < t_n=T$. 
We use equally spaced $t_i$ in our tests, but one could consider
non-uniform choices. 
If $U_t$ is a Brownian motion, then the increments 
\be
X_j=U_{t_j}-U_{t_{j-1}}
\ee
are independent and each is normal with mean zero and variance
$\kappa (t_j-t_{j-1})$. We test this joint distribution with a $\chi^2$ 
goodness-of-fit test. The idea is to divide the possible values of 
$(X_1,X_2,\cdots,X_n)$ into $m$ cells and count the 
number of samples that fall into each cell. 
Under the hypothesis that the $X_j$ are independent 
and normally distributed, we can compute the expected number of samples that 
fall into each cell. 
Let $O_j$ be the number of samples in cell $j$, and 
$E_j$ the expected value of this random variable under the 
hypothesis that $U_t$ is Brownian motion.
The test is then based on the statistic:
\be
    \chi^2 = \sum_{j=1}^m {(O_j-E_j)^2 \over E_j} 
\ee
Under the hypothesis that $U_t$ is a Brownian motion, 
if $m$ is large the distribution of this statistic is approximately the
$\chi^2$ distribution with $m-1$ degrees of freedom.

\begin{figure}[tbh]
\includegraphics{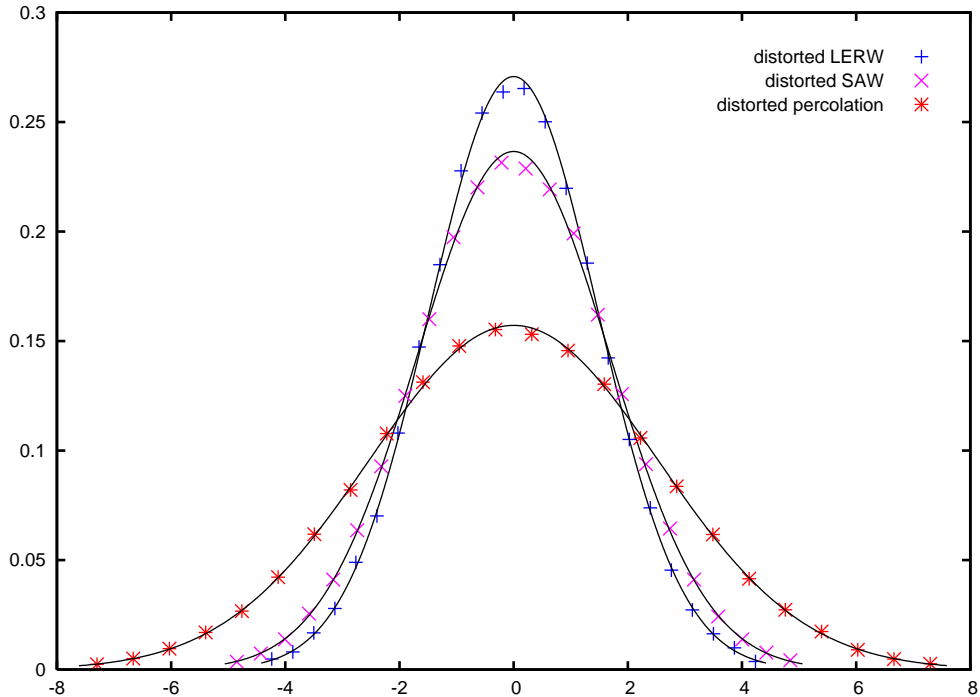}
\caption{\leftskip=25 pt \rightskip= 25 pt 
The points are a histogram for the density of $U_T/\sqrt{T}$ 
for the LERW, SAW and percolation with distortion $\lambda=0.95$. 
The curves are the density of a normal distribution with variance 
$\kappa$ where $\kappa$ is determined from the least squares fit.
}
\label{density}
\end{figure}

\begin{table}
\begin{center}
\begin{tabular}{|c|c|c|}
\hline
Model & $\lambda$ & $\kappa$ \\
\hline
LERW &   0.90 & 2.3315 $\pm$ 0.0106  \\
LERW &   0.95 & 2.1709 $\pm$ 0.0094  \\
LERW &   1.00 & 2.0008 $\pm$ 0.0093  \\
\hline
SAW &   0.90 & 3.0744 $\pm$ 0.0138  \\
SAW &   0.95 & 2.8414 $\pm$ 0.0108  \\
SAW &   1.00 & 2.6686 $\pm$ 0.0132  \\
\hline
percolation &   0.90 & 6.9446 $\pm$ 0.0287  \\
percolation &   0.95 & 6.4422 $\pm$ 0.0311  \\
percolation &   1.00 & 6.0404 $\pm$ 0.0265  \\
\hline
\end{tabular}
\caption{\leftskip=25 pt \rightskip= 25 pt 
The estimate of $\kappa$ using a weighted least squares fit for 
$E[U_t^2]$ as a function of $t$. The error bars are two standard deviations.
}
\label{kappa}
\end{center}
\end{table}

\begin{table}
\begin{center}
\begin{tabular}{|c|c|c|c|c|c|c|c|}
\hline
$\lambda$ & N &   $D(T/2)$ & $D(T)$ & $Z$  
& $\chi^2_a$ & $\chi^2_b$ & $\chi^2_c$ \\
\hline
 &  5,000 & 0.018695 & 0.223875 & 0.012146 & 0.306555 & 0.654153 & 0.139523\\
 & 10,000 & 0.605601 & 0.562645 & 0.012238 & 0.000614 & 0.000666 & 0.352671\\
0.90 
 & 20,000 & 0.578409 & 0.860558 & 0.029349 & 0.000000 & 0.000041 & 0.820988\\
 & 50,000 & 0.213470 & 0.308500 & 0.078267 & 0.000000 & 0.000000 & 0.000120\\
 & 100,000 & 0.222762 & 0.253184 & 0.003394 & 0.000000 & 0.000000 & 0.000000\\
\hline
\hline
 &  5,000 & 0.689088 & 0.567792 & 0.021239 & 0.913735 & 0.598610 & 0.766046 \\
 & 10,000 & 0.562020 & 0.898586 & 0.055965 & 0.423453 & 0.350098 & 0.593728 \\
0.95
 & 20,000 & 0.648585 & 0.638202 & 0.064346 & 0.793825 & 0.318862 & 0.865687 \\
 & 50,000 & 0.486077 & 0.178479 & 0.010384 & 0.078510 & 0.000173 & 0.850495 \\
 & 100,000 & 0.234619 & 0.004251 & 0.000812 & 0.000080 & 0.000000 & 0.501958 \\
\hline
\hline
 &  5,000 & 0.322835 & 0.909183 & 0.607075 & 0.599810 & 0.281570 & 0.332253 \\
 & 10,000 & 0.902220 & 0.499358 & 0.274931 & 0.118588 & 0.593618 & 0.820108 \\
1.00
 & 20,000 & 0.763553 & 0.856747 & 0.240621 & 0.167155 & 0.227972 & 0.322061 \\
 & 50,000 & 0.840997 & 0.746111 & 0.376101 & 0.695880 & 0.735182 & 0.702162 \\
 & 100,000 & 0.949877 & 0.934915 & 0.448990 & 0.305246 & 0.257668 & 0.509009 \\
\hline
\end{tabular}
\caption{\leftskip=25 pt \rightskip= 25 pt 
The distorted LERW. $\lambda$ is the amount of distortion, with $\lambda=1$ 
being no distortion. 
$N$ is the number of samples used. The other six columns give the p-value
of six different statistics used to test if the driving process is 
Brownian motion. See the text for details.
}
\label{tab_lerw}
\end{center}
\end{table}

\begin{table}
\begin{center}
\begin{tabular}{|c|c|c|c|c|c|c|c|}
\hline
$\lambda$ & N &   $D(T/2)$ & $D(T)$ & $Z$  
& $\chi^2_a$ & $\chi^2_b$ & $\chi^2_c$ \\
\hline

 &  5,000 & 0.076330 & 0.174654 & 0.538122 & 0.218454 & 0.170313 & 0.903646 \\
 & 10,000 & 0.471555 & 0.343790 & 0.824175 & 0.100552 & 0.019827 & 0.803437 \\
0.90
 & 20,000 & 0.246740 & 0.321011 & 0.905811 & 0.003206 & 0.000008 & 0.464170 \\
 & 50,000 & 0.165784 & 0.312940 & 0.013581 & 0.000000 & 0.000000 & 0.004110 \\
 & 100,000 & 0.251134 & 0.200176 & 0.000124 & 0.000000 & 0.000000 & 0.000000 \\
\hline
\hline
 &  5,000 & 0.408348 & 0.350261 & 0.933527 & 0.027830 & 0.037820 & 0.883070 \\
 & 10,000 & 0.136263 & 0.888704 & 0.681952 & 0.775130 & 0.704539 & 0.931479 \\
0.95
 & 20,000 & 0.805490 & 0.663008 & 0.983501 & 0.682962 & 0.797951 & 0.104660 \\
 & 50,000 & 0.483299 & 0.323554 & 0.882502 & 0.162307 & 0.030086 & 0.959757 \\
 & 100,000 & 0.313017 & 0.158268 & 0.794145 & 0.000018 & 0.009703 & 0.534966 \\
\hline
\hline
 &  5,000 & 0.456948 & 0.371554 & 0.636449 & 0.618423 & 0.126116 & 0.424437 \\
 & 10,000 & 0.557806 & 0.267715 & 0.821561 & 0.200096 & 0.145614 & 0.609878 \\
1.00
 & 20,000 & 0.933089 & 0.552806 & 0.514615 & 0.793825 & 0.253454 & 0.370379 \\
 & 50,000 & 0.956776 & 0.852597 & 0.588566 & 0.685205 & 0.275748 & 0.161052 \\
 & 100,000 & 0.501460 & 0.870474 & 0.219944 & 0.028353 & 0.079678 & 0.579479 \\
\hline
\end{tabular} \\
\caption{\leftskip=25 pt \rightskip= 25 pt 
The p-values of the distorted SAW.}
\label{tab_saw_1.00}
\end{center}
\end{table}

\begin{table}
\begin{center}
\begin{tabular}{|c|c|c|c|c|c|c|c|}
\hline
$\lambda$ & N &   $D(T/2)$ & $D(T)$ & $Z$  
& $\chi^2_a$ & $\chi^2_b$ & $\chi^2_c$ \\
\hline
&   5,000 &  0.268811 & 0.707407 & 0.035275 & 0.463308 & 0.440588 & 0.068530 \\
&  10,000 &  0.568890 & 0.306491 & 0.043429 & 0.387697 & 0.125325 & 0.286280 \\
0.90
&  20,000 &  0.113505 & 0.640373 & 0.192064 & 0.002475 & 0.000123 & 0.387076 \\
&  50,000 &  0.073017 & 0.749263 & 0.000274 & 0.000000 & 0.000000 & 0.007620 \\
& 100,000 &  0.000694 & 0.015332 & 0.000005 & 0.000000 & 0.000000 & 0.000000 \\
\hline
\hline
&   5,000 &  0.667251 & 0.576948 & 0.881721 & 0.978563 & 0.820268 & 0.261246 \\
&  10,000 &  0.831016 & 0.921401 & 0.964028 & 0.586700 & 0.601949 & 0.198512 \\
0.95
&  20,000 &  0.827859 & 0.938272 & 0.703387 & 0.580099 & 0.123268 & 0.124509 \\
&  50,000 &  0.504299 & 0.913135 & 0.465606 & 0.033206 & 0.000071 & 0.291786 \\
& 100,000 &  0.084361 & 0.515339 & 0.259067 & 0.000000 & 0.000000 & 0.043819 \\
\hline
\hline
&   5,000 &  0.075710 & 0.685162 & 0.242412 & 0.884633 & 0.296115 & 0.657406 \\
&  10,000 &  0.882136 & 0.304884 & 0.618193 & 0.267297 & 0.246049 & 0.766046 \\
1.00
&  20,000 &  0.271328 & 0.037469 & 0.133359 & 0.088622 & 0.587819 & 0.060798 \\
&  50,000 &  0.100113 & 0.008876 & 0.184139 & 0.256426 & 0.996808 & 0.283681 \\
& 100,000 &  0.333954 & 0.059206 & 0.518117 & 0.287868 & 0.779478 & 0.149034 \\
\hline
\end{tabular} 
\caption{\leftskip=25 pt \rightskip= 25 pt 
The p-values of the distorted percolation interface. }
\label{tab_perc_1.00}
\end{center}
\end{table}

We consider three choices of the cells. In the first choice we let $n=10$, 
and use only the signs of the $X_j$ to define the cells. Thus there 
are $2^{10}=1024$ cells. We denote the statistic in this case by 
$\chi^2_a$. 
The second choice of cells is similar - we use
only the signs of the increments but with $n=5$. So there are 
$2^5=32$ cells. This statistic is denoted by $\chi^2_b$. 
The third choice uses just two increments.
For each of the increments $X_1$ and $X_2$ we look at which quartile 
it falls in. We let $q$ be the constant such that under the hypothesis 
that $U_t$ is Brownian motion, 
$P(X_i>q)$=1/4, and then divide the possible values of each $X_i$ into 
the intervals $(-\infty,-q],[-q,0],[0,q],[q,\infty)$.
Thus there are $16$ cells. 
This final statistic is denoted by $\chi^2_c$. 
Note that the statistics $\chi^2_a$ and $\chi^2_b$ have the advantage 
that they do not involve the value of $\kappa$. For $\chi^2_c$ we need
a value for $\kappa$ to compute $q$. 

We study three values of the distortion parameter, $\lambda=0.9,0.95$
and $1$, for each of the three lattice models. 
Thus there are a total of nine models considered. 
(We have run the same simulations for $\lambda=1.05$ and $1.1$, but
do not present these results. They are consistent with the results that 
we do present.) 
We generate $10^5$ samples for each of the nine cases, 
compute the driving function $U_t$ for each sample, and 
record the value of $U_t$ at ten equally spaced times, 
$t=T/10,2T/10,\cdots,T$. 
We compute the variance of $U_t$ at each of the ten times
and then do a least squares fit to estimate $\kappa$, the slope
of $t \rightarrow E[U_t^2]$. The results are shown in table \ref{kappa}.

For each of the three models we present the results of our 
statistical tests in a table. 
We perform the statistical tests for the first $N$ samples of the 
$100K$ samples, using $N=5K, 10K, 20K, 50K$ and $100K$. 
In the tables we do not give the value
of the statistic, but rather the corresponding ``p-value.'' 
The p-value is defined as follows. Consider a one-sided statistic such as 
the Kolmogorov-Smirnov $D$ statistic. Suppose that in our test 
the value of the statistic is $x$. Then the p-value is the probability
$P(D>x)$ under the null hypothesis. For a two-sided statistic such as $Z$
the definition is modified in the obvious way. 
A small p-value (less than a few percent) indicates
that the value of the corresponding statistic is very unlikely under
the hypothesis that $U_t$ is Brownian motion, and so 
we should reject the hypothesis that $U_t$ is Brownian motion.

These nine tables contain a lot of numbers, but they clearly show the 
following. For $\lambda=1$ the p-values in the table for all three 
lattice models are not small,
and so our tests do not indicate that we should reject the hypothesis 
that $U_t$ is a Brownian motion. This is as it should be. 
When $\lambda=1$ the models should have a scaling limit
given by $SLE$, and so $U_t$ should indeed be Brownian motion. 

For the models that are not SLE, $\lambda=0.9,0.95$, we first consider 
the two Kolmogorov-Smirnov tests that $U_t$ is normal. The corresponding 
p-values are not typically small, and this test is completely 
ineffective at indicating that $U_t$ is not Brownian motion, even when 
we use $10^5$ samples. In figure \ref{density} the data points are  
histograms for the density function of $U_T/\sqrt{T}$ 
for the LERW, SAW and percolation with $\lambda=0.95$. The curves 
are the density functions for the normal 
distribution with variance $\kappa$ where $\kappa$ is determined
from our least squares fit. As the Kolmogorov-Smirnov test showed, 
the data points are fit very well by the normal curves. 
(For $\lambda \ne 1$, it is easy to show that the distorted model
is not SLE, and so $U_t$ is not a Brownian motion. However, this 
does not rule out the possibility that the $U_t$ are normal even 
for the distorted models.)

The other four tests involve the independence of the increments. 
The test based on the mean of the product of two independent increments
sometimes indicates correctly that $U_t$ is not Brownian motion, but
it is not very powerful. By contrast $\chi^2_a$ and $\chi^2_b$  
are quite effective at correctly indicating when $U_t$ is not a 
Brownian motion. For $\lambda=0.9$ these tests correctly indicate 
$U_t$ is not a Brownian motion with sample sizes on the order of ten to 
twenty thousand. For $\lambda=0.95$ these tests need on the order of 
a hundred thousand samples, but they are the only tests to correctly indicate 
that $U_t$ is not a Brownian motion for this amount of distortion.
The final statistic $\chi^2_c$ is only sometimes effective. 

\section{A faster zipper}
\label{zipper}

We briefly describe the standard method for computing the 
driving function of a simple curve $\gamma$. 
Let $g_s$ be the conformal map which takes the half plane minus 
$\gamma[0,s]$ onto the half plane, normalized so that for large $z$ 
\be
g_s(z) = z + { 2 t \over z} + O({1 \over z^2}),
\label{laurent_norm}
\ee
The coefficient $2t$ depends on $s$ and is 
the half-plane capacity of $\gamma[0,s]$.
The value of the driving function at $t$ is $U_t=g_s(\gamma(s))$.
Thus computing the driving function essentially reduces to 
computing this uniformizing conformal map. We will describe the 
``zipper algorithm'' for doing this \cite{kuh,mr}. 
Another approach to computing the driving function may be found in 
\cite{tsai}.

We find it more convenient to work with the conformal map
\be
h_s(z)=g_s(z)-U_s
\ee
It maps $\half \setminus \gamma[0,s]$ onto $\half$ and sends the tip
$\gamma(s)$ to the origin. The value of the driving function at $s$
is minus the constant term in the Laurent expansion of 
$h_s$ about $\infty$. From now on we work with this 
normalization for our conformal maps.

Let $z_0,z_1,\cdots,z_n$ be points along the curve with $z_0=0$. 
In our applications these are lattice sites. 
The zipper algorithm finds a sequence of 
conformal maps $h_i$, $i=1,2,\cdots,n$ such that 
$h_k \compose h_{k-1} \compose \cdots \compose h_1$
approximates the conformal map for the curve up to site $z_k$. 
Suppose that the conformal maps $h_1,h_2,\cdots,h_k$  have been 
defined so that $h_k \compose h_{k-1} \compose \cdots \compose h_1$
sends $\half \setminus \gamma$ to $\half$ where 
$\gamma$ is some curve that passes through $z_0,z_1, \cdots z_k$. 
In particular $z_k$ is mapped to the origin.
Let 
\be
w_{k+1}=h_k \compose h_{k-1} \compose \cdots \compose h_1(z_{k+1})
\label{compose}
\ee
Then $w_{k+1}$ is close to the origin.
We define $h_{k+1}$ to be a conformal map with the appropriate 
normalizations that sends $\half \setminus \gamma_{k+1}$ to $\half$
where $\gamma_{k+1}$ is a short simple curve that ends at $w_{k+1}$. 
The key idea is to choose this curve so that $h_{k+1}$
is explicitly known. The two choices we will use are 
``tilted slits'' and ``vertical slits.'' 

Let $2 \Delta t_i$ be the capacity of the map $h_i$, and $\Delta U_i$ 
the final value of the driving function for $h_i$. So 
\be
h_i(z) = z - \Delta U_i + { 2 \Delta t_i \over z} + O({1 \over z^2})
\ee
Then 
\be
h_k \compose h_{k-1} \compose \cdots \compose h_1(z) =
z  - U_t + { 2 t \over z} + O({1 \over z^2})
\ee
where 
\be
t= \sum_{i=1}^k \Delta t_i 
\ee
\be
U_t= \sum_{i=1}^k \Delta U_i 
\ee
Thus the driving function of the curve is obtained by ``adding up'' the 
driving functions of the elementary conformal maps $h_i$. 

We now consider the two particular types of maps we use for $h_{k+1}$. 
For tilted slits, $\gamma_{k+1}$ is the line segment from the origin
to $w_{k+1}$. There is no explicit formula for $h_{k+1}$ in this case, 
but there is a formula for its inverse:
\be
h_{k+1}^{-1}(z)=(z+x_l)^{1-\alpha} (z-x_r)^\alpha
\ee
where $x_l,x_r>0$. 
It maps the half plane onto the half plane minus a line segment which 
starts at the origin and forms an angle $\alpha$ with the positive real 
axis. The interval $[-x_l,x_r]$ gets mapped onto the slit. 
We must choose $x_l$ and $x_r$ so 
that $h_{k+1}$ satisfies our normalization conditions.
In particular, $h_{k+1}^{-1}$ must send the origin to the tip of the 
line segment, i.e., $w_{k+1}$. 
Tedious but straightforward calculation shows if we let 
$w_{k+1}=r\exp(i \alpha \pi)$, then 
\be
x_l=r \left({1-\alpha \over \alpha}\right)^\alpha, \qquad
x_r=r \left({\alpha \over 1-\alpha}\right)^{1-\alpha}
\ee
The changes in the driving function are given by 
\be
\Delta t={1 \over 4} r^2 \alpha^{1-2\alpha} (1-\alpha)^{2\alpha-1}, \qquad
\Delta U=r (1-2\alpha) \alpha^{-\alpha} (1-\alpha)^{-(1-\alpha)}
\ee

For vertical slits we take $\gamma_{k+1}$ to be the vertical line 
segment from the real axis to $w_{k+1}$. 
Note that since this $\gamma_{k+1}$ does not start at the origin,
this method does not approximate the original simple curve $\gamma$ by 
another simple curve. Instead, the domain of the conformal map that 
we construct to approximate $h_s$ is of the form $\half$ minus a set that 
is more complicated than a simple curve. 
This may make the reader nervous, but we will see in the next 
section that we get essentially the same driving functions 
using tilted slits or vertical slits.
The conformal map that removes this vertical line with our normalizations
is 
\be
h_{k+1}(z)=i\sqrt{-(z-x)^2-y^2}
\ee
where $w_{k+1}=x+iy$ and the branch cut for the square root is the 
negative real axis. 
The changes in the driving function are given by 
\be
\Delta t= {1 \over 4} y^2, \quad \Delta U=x
\ee

A comment on terminology is in order. We use ``zipper algorithm'' 
to refer to all the various algorithms we can get from different 
choices of the curve $\gamma_{k+1}$. Marshall and Rohde \cite{mr}
use ``zipper'' to refer only to the choice using tilted slits. 

The number of operations needed to compute a single $w_{k+1}$ is
proportional to $k$. So to compute all the points $w_{k+1}$
requires a time $O(N^2)$. To do better we must avoid evaluating 
the $k$-fold composition in \reff{compose} every time we compute a $w_{k+1}$. 
We begin by grouping the functions in \reff{compose} into
blocks. We denote the number of functions in a block by $\blength$. 
Let
\be
H_j = 
h_{jb} \compose h_{jb-1} \compose \cdots \compose h_{(j-1)b+2} \compose 
h_{(j-1)b+1} 
\label{blockdef}
\ee
If we write $k$ as $k=mb+r$ with $0 \le r < b$, then \reff{compose}
becomes
\be
w_{k+1} = 
h_{mb+r}  \compose h_{mb+r-1} \compose \cdots \compose h_{mb+1} \compose 
H_m \compose H_{m-1} \compose \cdots \compose H_1 (z_{k+1})
\label{blockcompose}
\ee
Typically, the number of compositions in \reff{blockcompose} 
is smaller than the 
number in \reff{compose} by roughly a factor of $b$. 
The $h_i$ are relatively simple, but the composition $H_j$ 
cannot be explicitly computed. 
Our strategy is to approximate the $h_i$ by functions 
whose compositions can be explicitly computed to give an explicit
approximation to $H_j$. This allows us to compute the compositions 
in \reff{blockdef} just once rather than every time we compute a $w_k$. 

Recall that $h_i$ is normalized so that $h_i(\infty)=\infty$ and 
$h_i^\prime(\infty)=1$. It maps $\half$ minus a simple curve which starts 
at the origin to $\half$. Let $h$ denote such a conformal map. 
Let $r$ be the largest distance from the origin to a point on the curve.
Then $h$ is analytic on $\{z \in \half: |z|>r\}$. 
Note that $h$ is real valued on the real axis.
By the Schwarz reflection principle it may be analytically continued 
to $\{z \in \complex: |z|>r\}$. Moreover, it does not vanish on this domain. 
So if we let $f(z)=h(1/z)$, then $f$ is analytic in 
$\{z \in \complex: |z|<1/r\}$ and $f(0)=0$, $f^\prime(0)=1$. 
The Laurent series of $h$ about $\infty$ is just the power series 
of $f$ about $0$. For large $z$, $h(z)$ is well approximated 
by a finite number of terms in this Laurent series. 
It will prove more convenient to work with a different series. 

Define $\hat{h}(z)=1/h(1/z)$. Since $h(z)$ does not vanish on 
$\{ z \in \complex: |z|>r\}$, $\hat{h}(z)$ is analytic in 
$\{ z \in \complex: |z| < 1/r \}$.
Our assumptions on $h$ imply that $\hat{h}(0)=0$ and $\hat{h}^\prime(0)=1$.
So $\hat{h}$ has a power series of the form 
\be
\hat{h}(z) = \sum_{j=1}^\infty \, a_j z^j
\label{hps}
\ee
with $a_1=1$. 
The radius of convergence of this power series is $1/r$.
Note that the coefficients of this power series are the
coefficients of the Laurent series of $1/h$. 

The primary advantage of working with the power series of $\hat{h}$
is its behavior with respect to composition. It is trivial to check that 
\be
(h_1 \compose h_2) \, \hat{} \, = \hat{h_1} \compose \hat{h_2}
\label{composeprop}
\ee
Our approximation for $h_i(z)$ is to replace $\hat{h_i}(z)$ 
by the truncation of its power series at order $\nterms$. So 
\be
h_i(z) = {1 \over \hat{h_i}(1/z)}
\approx \left[ \sum_{j=1}^n \, a_j z^{-j} \right]^{-1}
\ee

For each $h_i$ we compute the power series of $\hat{h_i}$ to order $\nterms$. 
We use them and \reff{composeprop} to compute the 
power series of $\hat{H_j}$ to 
order $\nterms$. Let $1/R_j$ be the radius of convergence for the
power series of $\hat{H_j}$. 
Now consider equation \reff{blockcompose}. 
If $z$ is large compared to $R_j$, then $H_j(z)$ is well approximated 
using the power series of $\hat{H}_j$.
We introduce a parameter $L>1$ and 
use this series to compute $H_j(z)$ whenever 
$|z| \ge L R_j$. When $|z| < L R_j$, we just use \reff{blockdef} to 
compute $H_j(z)$. The argument of $H_j$ is the result of applying the 
previous conformal maps to some $z_{k+1}$, and so is random. Thus whether 
or not we can approximate a particular $H_j$ by its 
series depends on the randomness and on which $w_{k+1}$ we are computing.

We need to compute $R_j$. Consider 
the images of $z_{(j-1)b},z_{(j-1)b+1}, \cdots z_{jb-1}$ 
under the map $H_{j-1} \compose H_{j-2} \compose \cdots \compose H_1$. 
The domain of the conformal map $H_j$ 
is the half-plane $\half$ minus some curve $\Gamma_j$ which passes through
the images of these points. The radius $R_j$ 
is the maximal distance from the origin to a point on $\Gamma_j$. 
This distance should be very close to or even equal to the maximum
distance from the origin to images of 
$z_{(j-1)b},z_{(j-1)b+1}, \cdots z_{jb-1}$ under 
$H_{j-1} \compose H_{j-2} \compose \cdots \compose H_1$. 
So we take $R_j$ to be the maximum of these distances.

Our algorithm depends on three parameters. The integer $\blength$ is 
the number of functions in a block. The integer $\nterms$ is the order 
at which we truncate the power series of the $\hat{H}_j$. 
The real number $L>1$ determines when we use the power series 
approximation for the block function. 

\begin{figure}[tbh]
\includegraphics{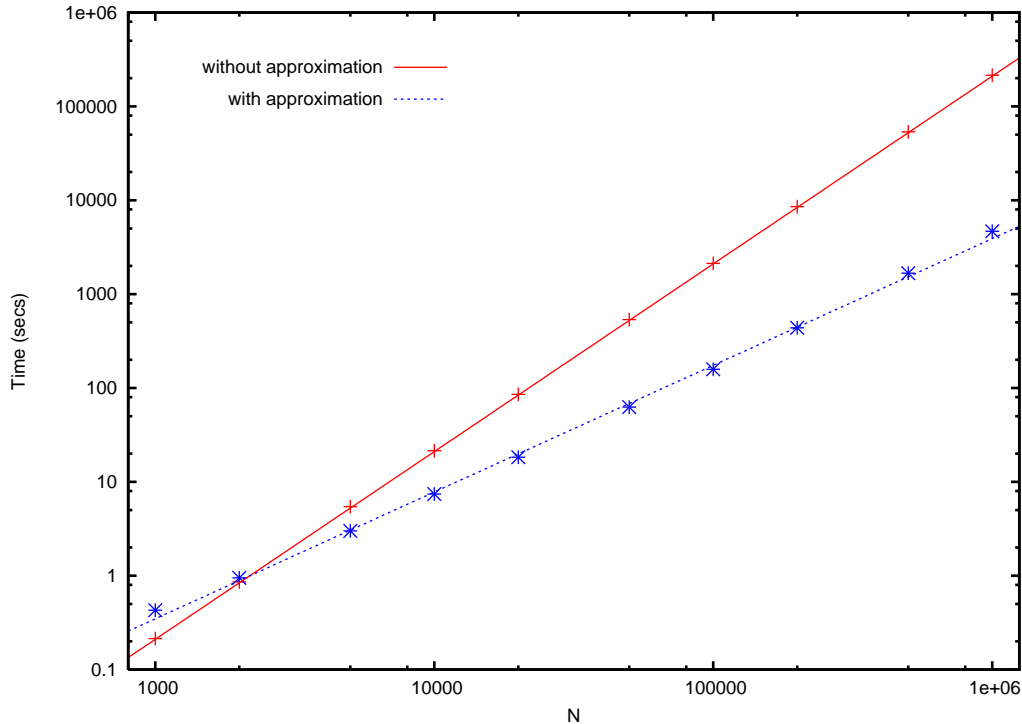}
\caption{\leftskip=25 pt \rightskip= 25 pt 
The points are the time (in seconds) needed to unzip a SAW with $N$ 
steps with and without the power series approximation. The lines have 
slopes $2$ and $1.35$. 
}
\label{loewner_timing}
\end{figure}

The improvement in the speed of the zipper algorithm from using our
power series approximation is shown in table \ref{table_laurent_timing} 
and figure \ref{loewner_timing}.
In these timing tests we use a single SAW with one million steps. 
We time how long it takes to unzip the first $N$ steps with 
and without the power series approximation. 
We do the computations using the power series approximation 
for different choices for the block length, namely 
$b=20,30,40,50,75,100,200,300$, and report the 
fastest time. The last column in the table indicates the block length
that achieves the fastest time. As a rule of thumb, a good choice 
for the block length (at least for the SAW) is $b=\sqrt{N}/4.$
The next to last column in the table gives the factor by which the 
use of the power series approximation reduces the time needed for 
the computation.
These timing tests were done on a PC with a 3.4 GHz Pentium 4 processor.

Without the power series approximation the time is $O(N^2)$. 
This is seen clearly in the log-log plot in figure \ref{loewner_timing} 
where the data for unzipping without the power series approximation 
is fit quite well by a line with slope 2. 
The data for unzipping using the power series approximation is 
fit by a line with slope 1.35. This indicates that the time required 
when the power series are used is approximately $O(N^{1.35})$.  

\begin{table}
\begin{center}
\begin{tabular}{|r|r|r|r|r|}
\hline
N & time 1 & time 2 & factor & block length  \\
\hline
    1,000 &  0.21 &  0.43 &  0.50  &  20 \\ 
    2,000 &  0.86 &  0.95 &  0.91  &  20 \\ 
    5,000 &  5.44 &  3.00 &  1.81  &  20 \\ 
   10,000 & 21.44 &  7.41 &  2.89  &  30 \\ 
   20,000 & 85.65 & 18.31 &  4.68  &  40 \\ 
   50,000 & 534.8 &  62.6 &  8.54  &  50 \\ 
  100,000 &  2128 &   158 & 13.45  &  75 \\ 
  200,000 &  8562 &   437 & 19.59  &  100 \\ 
  500,000 & 53516 &  1674 & 31.98  &  200 \\ 
1,000,000 & 214451 &  4675 & 45.87 &  200 \\ 
\hline
\end{tabular}
\caption{\leftskip=25 pt \rightskip= 25 pt 
The time (in seconds) needed to unzip a SAW with $N$ steps without 
using the power series approximation is shown
in the second column (time 1)
The time using the power series approximation is shown
in the third column (time 2).
The fourth column (factor) is the ratio of these two times. The block length
used is in the last column.
}
\label{table_laurent_timing} 
\end{center}
\end{table}

\section{Comparisons of computational methods}
\label{compare}

Given a simple curve, there are several choices when we compute 
its Loewner driving function. For the explicit conformal map in the zipper
algorithm we could use either the map that removes a vertical slit from the 
half plane or the map that removes a tilted slit. We could use the power
series approximation to speed up the calculation as explained in 
the last section.
We have a choice of how many points we use along the curve. 
In this section we study the effect of these various choices. 

For our study we generate a collection of forty self-avoiding walks with 
one million steps. For each walk we compute its driving function by three
different methods. The first uses the vertical slit map and the power
series approximation. The second uses the tilted slit map and the 
power series approximation. The third uses the vertical slit map 
without the power series approximation. To study the effect of using 
different numbers of points along the curve we do these computations
on subsets of the million points on the SAW. 
The subsets are obtained by taking every $n$th point along the walk.  
For the vertical slit map with the power series approximation 
(the fastest of the three methods) we use 
$n=1,2,5,10,20,50,100$. So the curves being unzipped have from $1,000,000$ 
to $10,000$ points on them. For the tilted slit map with the power series
approximation we use $n=2,5,10,20,50,100$. 
(This method is slower than that using vertical
slits because of the need to use Newton's method in the computation of 
the tilted slit conformal map.) For the vertical slit map without the 
power series approximation (by far the slowest of the three methods) 
we use $n=5,10,20,50,100$. 

Our choice of how to compare the driving functions computed by different
methods for the same curve merits some discussion. One might measure 
the difference between two driving functions by computing the 
supremum norm or the $L^1$ norm of their difference over a bounded interval. 
We do not do this.
The driving functions we are computing are approximations to Brownian 
motion sample paths. In particular, their slopes can be quite large. 
If you translate such a function by a small amount, the difference 
between the translated function and the original function can have 
a supremum norm or $L^1$ norm that is rather large. (Of course, the driving 
functions are continuous so these norms of the difference go to 
zero, but not linearly with the size of the translation.)
Thus a small error in computing the capacity produces a relatively
large error in these norms.
Instead of using these norms, we compare driving functions by using 
only the last point on the driving function. 
We denote this last point by $(T,U_T)$. The capacity of the SAW is 
$2T$ and $U_T$ is the image of the end of the SAW under the conformal
map that uniformizes the half plane minus the SAW.

We do not know the exact driving function of the SAW, so we 
treat the result of our computation using all one million steps (i.e., 
$n=1$) with the vertical slit map and the power series approximation 
as the exact answer. We then compute relative errors for $T$ and 
$U_T$. For $T$ we obtain the relative error by dividing the error in $T$ by
$T$. For $U_T$ we divide the error by the maximum of 
$|U_t|$ along the curve. For both of these relative errors we take the 
average over the forty SAW's. 

\begin{figure}[tbh]
\includegraphics{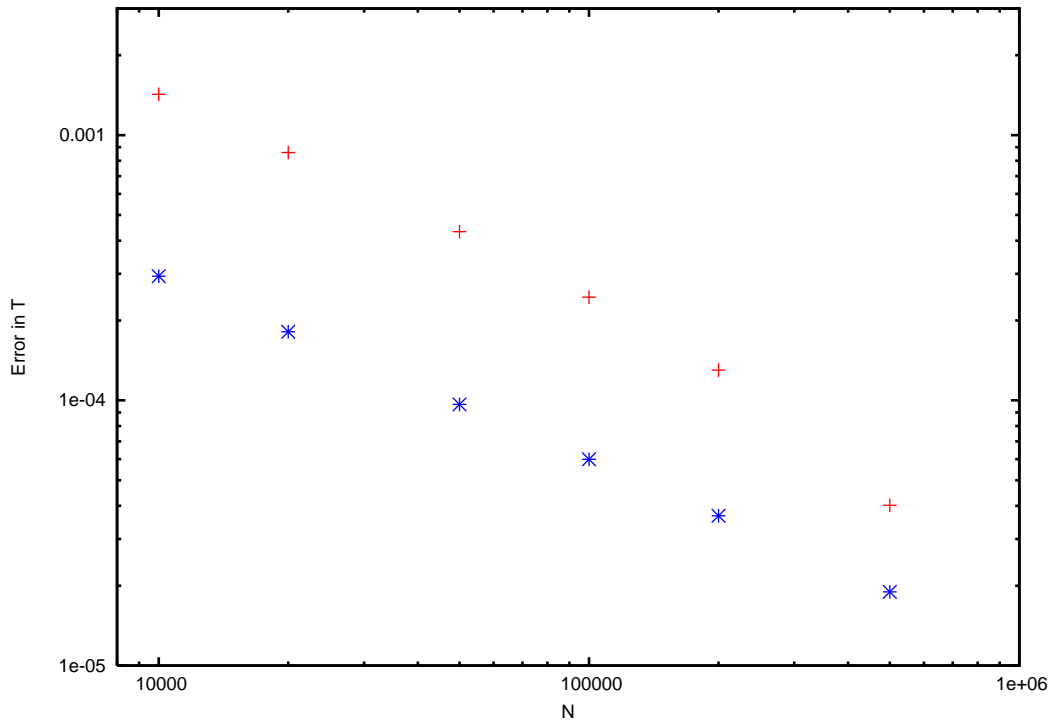}
\caption{\leftskip=25 pt \rightskip= 25 pt 
The top set of points is the relative error in $T$ as a function of 
the number of points $N$ used on the SAW. The bottom set of points
is the relative difference between $T$ computed using the tilted slit and 
vertical slit maps. 
}
\label{loewner_cap}
\end{figure}

\begin{figure}[tbh]
\includegraphics{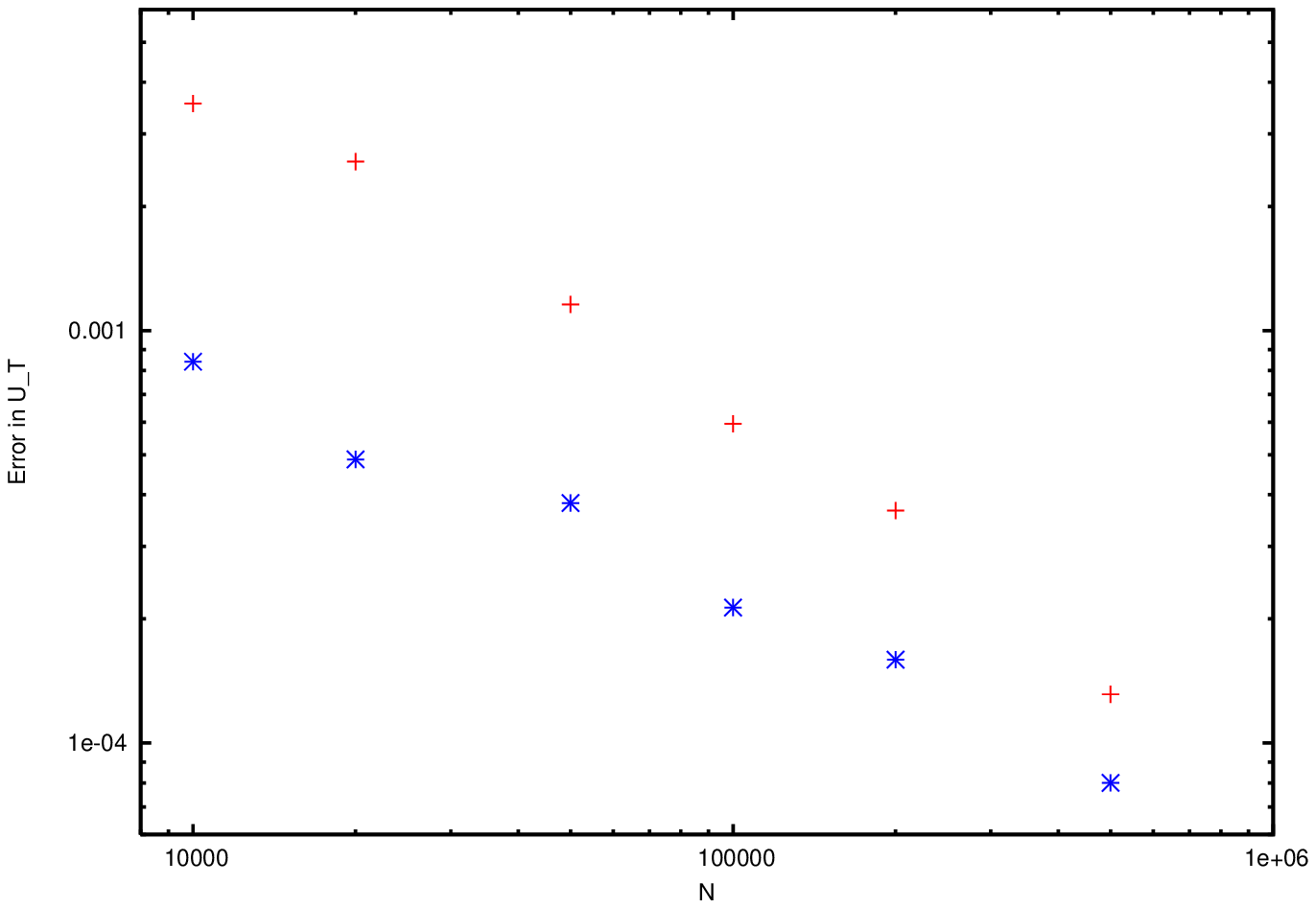}
\caption{\leftskip=25 pt \rightskip= 25 pt 
The top set of points is the relative error in $U_T$ as a function of 
the number of points $N$ used on the SAW. The bottom set of points
is the relative difference between $U_T$ computed using the tilted slit and 
vertical slit maps. 
}
\label{loewner_drive}
\end{figure}

We begin with the effect of the number of points we use along the 
curve. We study this using the vertical slit map with the power series
approximation. 
We let $N=10^6/n$ denote the number of points on the SAW used in
a particular computation. For the vertical slit map with the power
series approximation, $N$ ranges from $10^4$ to $10^6$. 
In figure \ref{loewner_cap} we plot the relative error in $T$
as a function of $N$. This is the higher of the two sets of points.
In figure \ref{loewner_drive} we plot the relative error in $U_T$
as a function of $N$. Again, this is the higher of the two sets of points.
In both of these plots the error is computed by treating the values for 
$N=10^6$ as the exact values.
These plots give an idea of the size of the 
error when the number of points used on the curve is reduced. 

To study the effect of the choice of conformal map (vertical slit vs. 
tilted slit), we compute the difference between the values of T 
using the two different maps with the same value of $N$. We do the same
for $U_T$. We convert these errors to relative errors by dividing by $T$ 
or the maximum of $|U_t|$. 
These relative errors are plotted as functions of 
$N$ in figures \ref{loewner_cap} and \ref{loewner_drive}.
In both figures they are the lower sets of points.
The figure shows they go to zero as $N \rightarrow \infty$. 
Note that these points lie well below the points that show the 
difference between the vertical map computation at the given $N$
with the vertical map computation with $N=10^6$. In other words, 
the effect of the choice of the conformal map is smaller than the 
effect of using more points along the curve. 

\begin{table}
\begin{center}
\begin{tabular}{|l|c|c|c|c|}
\hline
N & $T$ & $U_T$  \\
\hline
 10,000  & 3.66e-09  &   9.12e-09  \\
 20,000  & 2.59e-09  &   7.98e-09  \\
 50,000  & 1.61e-09  &   4.73e-09  \\
 100,000 & 8.38e-10  &   3.15e-09  \\ 
 200,000 & 6.98e-10  &   1.76e-09  \\ 
\hline
\end{tabular}
\caption{\leftskip=25 pt \rightskip= 25 pt 
The effect of the power series approximation. The table shows the 
relative differences in $T$ and in $U_T$ when we use 
and do not use the power series approximation.
}
\label{table_laurent} 
\end{center}
\end{table}

Finally we discuss the effect of the power series approximation. 
We use the vertical slit method. We compute the differences 
in $T$ and $U_T$ computed with and without the power series approximation. 
We convert these differences to relative differences by dividing by $T$ 
or the maximum of $|U_t|$. 
Table \ref{table_laurent} shows the relative differences. 
They are very small and insignificant compared to the 
differences that we see when we change the conformal map used 
or the number of points along the curve. In this study we took the number
of terms in our power series to be $12$ and the parameter $L=4.0$.
Increasing either of these improves the accuracy of the approximation
while slowing down the computation. We used block lengths of $100$. 
 
\section{Conclusions}

In this paper we computed the stochastic driving process of  
several models of random curves which we know are not SLE. We 
considered several statistical tests of whether this driving process
is a Brownian motion. Simply checking if the distribution at a fixed 
time is Gaussian was seen to be useless. We must use 
a test that involves the independence of the increments of the 
Brownian motion. Our most effective tests were 
$\chi^2$ goodness of fit tests in which we consider $n$ equal increments 
of the process and take the cells to be determined simply by the signs
of the increments. This test was the most successful at concluding 
that for the models that are not SLE, the stochastic driving process
is not a Brownian motion. One nice feature of this test is that it 
does not involve the value of $\kappa$. 

In the models we studied which are not an SLE, we have broken 
conformal invariance in a drastic way. In particular, these models
are not locally isotropic. The models from physics that have 
been recently studied as possible SLE's should 
be locally isotropic. So a lack of conformal invariance in these
physical models would have to arise in a completely different
way compared to the models we studied. 
Our main conclusion is not that the specific 
tests of Brownian motion that we found effective for our models 
are the best tests for all models, but rather that when we test
for SLE by computing the Loewner driving process, we should employ a 
variety of tests that this process is a Brownian motion.

We have also studied the numerical problem of computing the driving
function of a given curve using the zipper algorithm. 
We have seen that the difference in the driving function when we use
vertical slits or tilted slits for the elementary conformal map for the 
algorithm is quite small. 
Given that the vertical slit map is considerably faster and easier 
to implement, we see no reason to use the tilted slit map. 
We have also shown that the speed of this algorithm 
can be increased dramatically using power series approximations 
of certain analytic functions. The loss of accuracy from this approximation
is extremely small, in particular it is insignificant compared to the 
effect of changing the number of points used to define the curve we 
are unzipping or compared to the difference between using vertical slits
or tilted slits in the algorithm. 

\begin{appendix}
\section{Details of the simulations}

In this appendix we give some details of the simulations of the distorted 
LERW, SAW and percolation interface.

The LERW walk that we simulate is chordal LERW in the half plane from 
$0$ to $\infty$. This means that we take an ordinary random walk 
beginning at the origin and condition it to remain in the upper half 
plane. Then we erase the loops in chronological order. The ordinary 
random walk conditioned to remain in the upper half plane is 
easy to simulate since it is 
given by a random walk beginning at $0$ with transition probabilities 
that only depend on the vertical component of the present location 
of the walk. If the site has vertical component $k$, then the walk 
moves up with probability $(k+1)/4k$, down with probability $(k-1)/4k$, 
and to the right or left with probability $1/4$. 
(See, for example, section 0.1 of \cite{lawler}).
This process is known as the half plane excursion. 
The half-plane excursion is transient, i.e., each lattice site is visited 
by the excursion a finite number of times. This implies that the loop erasure
makes sense. (For a transient walk all parts of the walk would eventually
be part of a loop and so would be erased.)
Note, however, that if we take an infinite half plane excursion and only 
consider the first $n$ steps and loop-erase this walk, the result does not 
completely agree with the loop-erasure of the full infinite excursion.
A site which is visited by the excursion before time $n$ may be erased 
by a loop formed after time $n$. 

In practice there is no way to know if a visit to a site will be erased 
by some future loop without simulating the entire excursion. 
So in the simulation we do the following. 
We generate a half-plane excursion, erasing the loops as 
they are formed. We stop when the resulting walk has $N$ steps. 
If $n$ is small compared to $N$, 
then the distribution of our walk for the first $n$ steps is 
close to the true distribution of the first $n$ steps of the LERW. 
We will only compute the driving function for the first $n$ steps.

We take $N=50,000$ and generate $100,000$ samples.  
If we work on a unit lattice, a LERW with $N$ steps 
has a size of order $N^\nu$ with $\nu=4/5$. So we rescale our walk 
by a factor of $N^\nu$ to obtain a curve whose size is of order one. 
We then compute its driving function up to time $T=0.01$. 
The time $T$ is one half of the capacity. So the number of steps
needed to reach $T=0.01$ is random. 
For this choice of $T$, the mean of this random number of steps is 
approximately $8200$, roughly a factor of six smaller than $N$. 

The SAW in the upper half plane is defined as follows. Let $N$ be a positive
integer. We consider all nearest neighbor walks with $N$ steps
in the upper half plane which begin at the origin and do not visit
any site more than once. We put the uniform probability measure on this
finite set of walks. We let $N \rightarrow \infty$ to get a 
probability measure on infinite self-avoiding walks on the unit lattice
in the upper half plane. Then we take the lattice spacing to zero.  
We simulate the SAW in the half plane with a fixed number of steps 
with the pivot algorithm, a Markov Chain Monte Carlo method 
\cite{ms}. We use the fast implementation of this algorithm 
introduced in \cite{tk_pivot}. For the SAW there is an issue similar to
the LERW. The pivot algorithm produces the uniform distribution on 
the set of walks with $N$ steps. But this is not the distribution of 
the infinite SAW in the half plane restricted to walks of length $N$. 
As with the LERW, we address this problem by simulating walks with $N$ steps
but then computing the driving function for only the first $n$ steps 
where $n$ is much smaller than $N$. 

We simulate SAW's with $200,000$ steps. We sample the SAW from the pivot 
algorithm every $100,000$ time steps in the Markov chain. We run the chain
for $10^{10}$ iterations to produce $100,000$ samples. Unlike the other 
two models, these samples are not exactly independent, but the large time
interval between sampling makes the samples very close to 
independent. We rescale the SAW by a factor of $N^\nu$ with 
$\nu=3/4$ and then compute its driving function up to time $T=0.002$. 
The mean of the number of steps needed to reach $T=0.002$ 
is approximately $9350$, roughly a factor of $20$ smaller
than the total number of steps in the SAW. 

The percolation model we study is site percolation on the triangular lattice
in the upper half plane, but we describe it using the  
hexagonal lattice in the upper half plane. 
Each hexagon is colored white or black with probability $1/2$. 
The hexagons along the negative real axis are white and those along the 
positive real axis are black. This forces an interface which starts 
with the bond through the origin between the adjacent differently 
colored hexagons on the real axis. This interface is the unique
curve on the hexagonal lattice which begins at this bond and 
has all white hexagons along one side of the interface and all black 
ones along the other side.

Note that unlike the LERW or SAW there is no finite length effect for 
percolation interfaces.
If we generate interfaces with $n$ steps, they have exactly the same 
distribution as the first $n$ steps of interfaces of length $N$ 
where $N>n$. 
We generate $100,000$ samples of interfaces with $N=40,000$ steps. 
We rescale our walk by a factor of $N^\nu$ with $\nu=4/7$ 
and then compute the driving function up to time $T=0.1$. 
This corresponds to a mean number of steps of approximately $11,300$. 

We end with a comment on the time $T$ and our rescaling of the various
curves. For each model we have rescaled the curves by a factor of $N^\nu$. 
This is merely for convenience. We could have left the curves on 
a unit lattice and computed the driving function up to time
$T$ given by the above values times $N^{2 \nu}$. 
What is important is that the mean number of steps of the curves we 
are unzipping is large (so that we are close to the scaling limit)
but still significantly smaller than the total number of 
steps in the curve (for the SAW and LERW) so that we avoid the 
finite length effects discussed above. We have chosen the values
of $T$ so that in all three models 
the mean number of steps unzipped is on the order of $10,000$. 
Almost all of the time in these simulations is spent on computing 
the driving functions. Generating the random curves takes 
essentially no time by comparison. 

\end{appendix}

\bigskip

\noindent {\bf Acknowledgments:}
This research was inspired by talks and interactions during a visit 
to the Kavli Institute for Theoretical Physics in September, 2006.
I thank Don Marshall and Stephen Rohde for  
useful discussions about the zipper algorithm.
This research was supported in part by the National Science Foundation 
under grant DMS-0501168.


\bigskip
\bigskip

\begin{thebibliography}{99}

\newcommand \jtype{\it}

\def \prs    {{\jtype Proc. R. Soc. London}\ }
\def \jsp    {{\jtype J. Statist. Phys.}\ }
\def \pr     {{\jtype Phys. Rev.}\ }
\def \prb    {{\jtype Phys. Rev. B}\ }
\def \prl    {{\jtype Phys. Rev. Lett.}\ }
\def \cmp    {{\jtype Commun. Math. Phys.} \ }
\def \jpc    {{\jtype J. Phys.: Condens. Matter}\ }
\def \jpa    {{\jtype J. Phys. A: Math. Gen.}\ }
\def \jpcold {{\jtype J. Phys. C: Solid State Phys.}\ }
\def \pl     {{\jtype Phys. Lett.}\ }
\def \lmp    {{\jtype Lett. Math. Phys.}\ }
\def \npb    {{\jtype Nucl. Phys. B}\ }
\def \jmp    {{\jtype J. Math. Phys.}\ }
\def \jap    {{\jtype J. Appl. Phys.}\ }
\def \jpsj   {{\jtype J. Phys. Soc. Jpn.}\ }
\def \rmp    {{\jtype Rev. Math. Phys.}\ }
\def \epl    {{\jtype Europhys. Lett.}\ }
\def \aihp   {{\jtype Ann. Inst. H. Poincar\'e}\ }
\def \jfa    {{\jtype J. Funct. Anal.}\ }
\def \zpcm   {{\jtype Z. Phys. B}\ }
\def \sc     {{\jtype Science}\ }

\bibitem{ahhm}
C.~Amoruso, A.~K.~Hartman, M.~B.~Hastings, and M.~A.~Moore,
\newblock{Conformal invariance and SLE in two-dimensional Ising spin glasses}, 
\prl {\bf 97}, 267202 (2006).
Archived as {\tt arXiv:cond-mat/0601711}.

\smallskip

\bibitem{bbcfA}
D.~Bernard, G.~Boffetta, A.~Celani, and G.~Falkovich,
\newblock{Conformal invariance in two-dimensional turbulence}, 
\newblock{\em Nature Physics} {\bf 2}, 124 (2006).
Archived as {\tt arXiv:nlin.CD/0602017}.
\smallskip

\bibitem{bbcfB}
D.~Bernard, G.~Boffetta, A.~Celani, and G.~Falkovich,
\newblock{Inverse turbulent cascades and conformally invariant curves}. 
Archived as {\tt arXiv:nlin.CD/0609069}.
\smallskip

\bibitem{bdm}
D.~Bernard, P.~Le Doussal, and A.~A.~Middleton,
\newblock{Are domain walls in 2D spin glasses described by stochastic
Loewner evolutions?}, \prb {\bf 76}, 020403(R) (2007).
Archived as {\tt arXiv:cond-mat/0611433}.
\smallskip

\bibitem{cn}
F.~Camia and C.~M.~Newman,
\newblock{Critical percolation exploration path and SLE(6): 
a proof of convergence}, preprint. 
Archived as {\tt arXiv:math.PR/0604487}.
\smallskip

\bibitem{tk_sle}
T.~Kennedy,
\newblock{A fast algorithm for simulating the chordal Schramm-Loewner 
evolution}, 
\newblock{\em J. Statist. Phys.} {\bf 128}, 1125--1137 (2007).
Archived as {\tt arXiv:math.PR/0508002}.
\smallskip

\bibitem{tk_pivot}
T.~Kennedy,
\newblock{A faster implementation of the pivot algorithm for self-avoiding 
walks},
\newblock{\em J. Statist. Phys.} {\bf 106}, 407--429 (2002).
Archived as {\tt arXiv:cond-mat/0109308}
\smallskip

\bibitem{kuh}
 R.~K\"uhnau, 
\newblock{Numerische Realisierung konformer Abbildungen durch 
"Interpolation"},
\newblock{\em Z. Angew. Math. Mech.} {\bf 63}, 631-­637 (1983).

\bibitem{lawler} 
G.~Lawler, 
\newblock{\it Conformally Invariant Processes in the Plane},
\newblock{\it Mathematical Surveys and Monographs, vol. 114},
\newblock{American Mathematical Society, 2005}.

\bibitem{lsw_saw}
G.~Lawler, O.~Schramm, and W.Werner, 
\newblock{On the scaling limit of planar self-avoiding walk}, 
{\it Fractal Geometry and Applications: a Jubilee of Benoit Mandelbrot, 
Part 2}, 339--364, \newblock{\em Proc. Sympos. Pure Math. 72}, 
Amer. Math. Soc., Providence, RI, 2004.
Archived as {\tt arXiv:math.PR/0204277}.

\bibitem{lsw_lerw}
G.~Lawler, O.~Schramm, and W.Werner, 
\newblock{Conformal invariance of planar loop-erased random walks and 
uniform spanning trees}, 
\newblock{\em Ann. Probab.} {\bf 32}, 939--995 (2004).
Archived as {\tt arXiv:math.PR/0112234}.

\bibitem{ms} 
N.~Madras and G.~Slade, 
\newblock{\it The Self-Avoiding Walk}, 
\newblock{Birkh\"auser, Boston-Basel-Berlin, 1993}. 

\bibitem{mr} 
D.~E.~Marshall and S.~Rohde,
\newblock{Convergence of a variant of the Zipper algorithm 
for conformal mapping}, \newblock{\em SIAM J. Numer. Anal.} to appear.

\bibitem{schramm} 
O.~Schramm,
\newblock{Scaling limits of loop-erased random walks and uniform spanning 
trees},
{\it Israel J. Math. }  {\bf 118}, 221--288 (2000).
Archived as {\tt arXiv:math.PR/9904022}.

\bibitem{smirnov_perc}
S.~Smirnov,
Critical percolation in the plane,
{\it C. R. Acad. Sci. Paris S\'er. I Math.} {\bf 333}, 239 (2001).

\bibitem{smirnov_ising}
S.~Smirnov,
\newblock{Towards conformal invariance of 2D lattice models},
Proceedings of the International Congress of Mathematicians, 
Vol. II, 1421-1451, Eur. Math. Soc., Zurich, 2006.
Archived as {\tt arXiv:0708.0032v1 [math-ph]}.
\smallskip

\bibitem{ss}
O.~Schramm and S.~Sheffield, 
\newblock{Contour lines of the two-dimensional discrete Gaussian free field},
preprint. Archived as {\tt arXiv:math.PR/0605337}.
\smallskip

\bibitem{tsai}
J.~Tsai,
\newblock{The Loewner driving function of trajectory arcs of 
quadratic differentials}, preprint. Archived as 
{\tt arXiv:0704.1933v2 [math.CV]}.
\smallskip

\bibitem{zhan}
D.~Zhan,
\newblock{The scaling limits of planar LERW in finitely connected domains},
\newblock{\em Ann. Probab.} {\bf 36}, 467--529 (2008).
Archived as {\tt arXiv:math.PR/0610304}.
\smallskip

\end{thebibliography}
\end{document}